\documentclass[11pt,reqno]{amsart}
\usepackage[T1]{fontenc}
\usepackage[utf8]{inputenc}
\usepackage{lmodern}
\usepackage{microtype}
\usepackage[a4paper,margin=29mm]{geometry}
\usepackage{amsmath,amssymb,mathtools,mathrsfs}
\usepackage{enumitem,booktabs}
\usepackage{xcolor}
\usepackage[hidelinks]{hyperref}
\hypersetup{pdftitle={Depth two and fiber functors},
  pdfauthor={Mainak Ghosh and S\'ebastien Palcoux}}
\setlist[enumerate]{label=\textup{(\roman*)},leftmargin=*,itemsep=4pt}
\setlist[itemize]{leftmargin=*,itemsep=3pt}
\newtheorem{theorem}{Theorem}[section]
\newtheorem{proposition}[theorem]{Proposition}
\newtheorem{lemma}[theorem]{Lemma}
\newtheorem{corollary}[theorem]{Corollary}
\theoremstyle{definition}
\newtheorem{definition}[theorem]{Definition}
\newtheorem{example}[theorem]{Example}
\theoremstyle{remark}
\newtheorem{remark}[theorem]{Remark}
\newcommand{\C}{\mathcal C}
\newcommand{\M}{\mathcal M}
\newcommand{\one}{\mathbf 1}
\newcommand{\id}{\mathrm{id}}
\newcommand{\Hom}{\operatorname{Hom}}
\newcommand{\End}{\operatorname{End}}
\newcommand{\Nat}{\operatorname{Nat}}
\newcommand{\Vect}{\operatorname{Vec}}
\newcommand{\Rep}{\operatorname{Rep}}
\newcommand{\Corep}{\operatorname{Corep}}
\newcommand{\add}{\operatorname{add}}
\newcommand{\ev}{\operatorname{ev}}
\newcommand{\coev}{\operatorname{coev}}
\newcommand{\op}{\mathrm{op}}

\newcommand{\boxt}{\mathbin{\boxtimes}}
\newcommand{\depth}{\operatorname{depth}}
\newcommand{\eps}{\varepsilon}
\newcommand{\xto}[1]{\xrightarrow{#1}}
\allowdisplaybreaks[2]
\numberwithin{equation}{section}
\title[Depth two and fiber functors]{Depth two and fiber functors}
\author{Mainak Ghosh}
\address{Hetao Institute of Mathematics and Interdisciplinary Sciences,
Futian District, Shenzhen, China}
\email{ghosh.main@gmail.com}
\author{S\'ebastien Palcoux}
\address{Beijing Institute of Mathematical Sciences and Applications,
Huairou District, Beijing, China}
\email{sebastienpalcoux@gmail.com}
\date{November 26, 2026}

\begin{document}
\begin{abstract}
We generalize the depth-two Hopf reconstruction principle of the Ocneanu--Szyma\'nski theorem to finite tensor categories over arbitrary fields. Such a category admits a fiber functor if and only if it is tensor-generated by a connected algebra object of depth at most two. Every such algebra is automatically Frobenius and becomes the regular comodule algebra of a finite-dimensional Hopf algebra. The main argument derives the absorbing family for all objects from one depth-two identity, without assuming semisimplicity. For unitary fusion categories, the characterization gives an equivalent criterion for on-site symmetry realizations on tensor-product spin chains.
\end{abstract}
\maketitle
\markboth{MAINAK GHOSH AND S\'EBASTIEN PALCOUX}{DEPTH TWO AND FIBER FUNCTORS}

\makeatletter
\begingroup
\let\saved@tocline\@tocline
\def\@tocline#1#2#3#4#5#6#7{%
  \saved@tocline{#1}{#2}{#3}{#4}{\small#5}{#6}{#7}}
\tableofcontents
\endgroup
\makeatother
\clearpage

\section{Introduction}

Depth two connects the structure of a subfactor with finite-dimensional Hopf symmetry \cite{Szymanski,Longo}. We generalize this Hopf-reconstruction principle of the Ocneanu--Szyma\'nski theorem to finite tensor categories over arbitrary fields: a single connected tensor-generating algebra of depth at most two determines a fiber functor and a regular Hopf-comodule model. The generalized conclusion is the categorical reconstruction of the Hopf symmetry. Connectedness corresponds to irreducibility, and Proposition~\ref{prop:subfactor-depth} identifies our depth convention with the rooted principal-graph convention. The role of Frobenius algebra objects and their two-categorical factorizations is developed in \cite{Muger,GP}.

Throughout, $k$ is a field and all categories and functors are $k$-linear. A \emph{finite tensor category} is a finite abelian rigid monoidal category with exact tensor product, simple unit $\one$, and $\End(\one)=k$. A \emph{fiber functor} is an exact faithful strong monoidal functor to the category $\Vect_k$ of finite-dimensional vector spaces. An algebra object $(A,m,e)$ is \emph{connected} if $\Hom(\one,A)=k e\cong k$. It \emph{tensor-generates} $\C$ if every object of $\C$ is a subquotient of a finite direct sum of powers $A^{\otimes r}$, including $A^{\otimes0}=\one$. Section~\ref{sec:prelim} gives the detailed definitions.

\subsection{Depth and its depth-two form}
The depth of an algebra is read from an alternating tower. In the standard two-categorical factorization $A=Y\circ Y^\vee$ associated with $A$, its levels are
\begin{equation}\label{eq:tower-intro}
\one,\quad Y,\quad A,\quad A\circ Y,\quad
A^{\otimes2},\quad A^{\otimes2}\circ Y,\quad\ldots,
\end{equation}
starting at level zero. The even and odd objects belong to two different Hom categories. The \emph{depth} is the last level at which a new isomorphism class of indecomposable direct summand appears, with comparisons made only within the same sector; it is infinite when such first appearances are unbounded. Definition~\ref{def:all-depth} gives the precise intrinsic formulation, and Section~\ref{sec:depth} constructs the factorization. For a Frobenius algebra, the two connecting morphisms are two-sided adjoints. For the canonical algebra of a finite-index, finite-depth subfactor, Proposition~\ref{prop:subfactor-depth} identifies its depth with principal-graph depth rooted at the trivial even bimodule.

Write $\add(Y)$ for the direct summands of finite direct sums of copies of $Y$. The condition needed here has the short reformulation
\begin{equation}\label{eq:depth-intro}
\depth(A)\leq2
\quad\Longleftrightarrow\quad A\circ Y\in\add(Y).
\end{equation}
If $A$ is connected, $\End(Y)=k$, and this is equivalently
\begin{equation}\label{eq:connected-depth-intro}
                    A\circ Y\cong Y^{\oplus n}\qquad(n\geq1).
\end{equation}
These statements, including why level three controls all subsequent levels, are proved in Proposition~\ref{prop:depth-two-identity}. The isomorphism in \eqref{eq:connected-depth-intro} is taken in the odd Hom category, not just between underlying even objects.

\subsection{The characterization}
\begin{samepage}
\begin{theorem}[Depth two and fiber functors]\label{thm:main}
Let $\C$ be a finite tensor category over $k$. The following conditions are equivalent.
\begin{enumerate}
\item $\C$ admits a fiber functor.
\item $\C$ is tensor-generated by a connected algebra object of depth at most two.
\end{enumerate}
Every algebra in \textup{(ii)} admits a Frobenius structure. Moreover, a suitable tensor equivalence $\C\simeq\Corep(H)$ identifies it, as an algebra object, with the regular comodule algebra of a finite-dimensional Hopf algebra $H$.
\end{theorem}
\end{samepage}

Here $\Corep(H)$ denotes the category of finite-dimensional left $H$-comodules. The qualification ``at most'' accounts only for the trivial case: the unit algebra tensor-generates $\Vect_k$ and has depth one. If $\C\not\simeq\Vect_k$, every algebra in \textup{(ii)} has depth exactly two; see Corollary~\ref{cor:exact-depth}. The Frobenius specialization gives the characterization proposed in \cite[Speculation~13.16]{GP}, with depth defined through indecomposable direct summands and tensor generation understood for all objects. The theorem requires no Frobenius or separability assumption on the algebra. The proof also constructs the fiber functor explicitly and shows that the algebra is a projective generator of $\C$; the full structural statement is Corollary~\ref{cor:structure}.

\subsection{Historical context}
The operator-algebraic starting point is the reconstruction associated with Ocneanu. Szyma\'nski \cite{Szymanski} established the Hopf crossed-product description for irreducible depth-two subfactors, while Longo \cite{Longo} expressed the correspondence through categories of endomorphisms. David \cite{David} gave a detailed treatment of Ocneanu's paragroup construction and its relation to finite-dimensional Kac algebras. Das--Kodiyalam \cite{DK} supplied a planar-algebra proof of the Ocneanu--Szyma\'nski theorem, and Kodiyalam--Sunder \cite{KS} treated the planar algebra of a semisimple and cosemisimple Hopf algebra. Ghosh \cite{Ghosh} developed the pivotal two-categorical viewpoint on planar algebras. In the algebraic theory of extensions, Kadison--Nikshych \cite{KN} related depth-two Frobenius extensions to Hopf actions and Galois extensions.

The closest categorical antecedents are M\"uger's depth-two reconstruction in semisimple spherical irreducible Morita contexts \cite[\S6.3]{Muger} and the absorbing-monoid construction of M\"uger--Tuset \cite[Proposition~3.2]{MT}. The latter constructs a fiber functor from a connected monoid that absorbs every object as a module in a semisimple tensor category. Our main step derives that entire family from the single depth-two identity, using a projective generator to control extensions; see Remark~\ref{rem:absorption}. Direct-summand similarity is already central to algebraic depth \cite[\S1.1 and Definition~2.3]{KadisonDepth}; we retain both alternating sectors and their first-appearance numbering.

We use the finite-tensor-category framework of Etingof--Ostrik \cite{EO} and Etingof--Gelaki--Nikshych--Ostrik \cite{EGNO}. After the module-category collapse, the fiber functor is supplied by the standard module-category correspondence \cite[Proposition~7.3.3 and Example~7.4.6]{EGNO}, followed by reconstruction \cite[Theorem~5.3.12]{EGNO}; see also \cite[\S\S2.2--2.3]{Schauenburg} for arbitrary-field reconstruction. For the regular Hopf-comodule algebra we use Ghosh--Palcoux \cite[Theorem~11.11 and Lemma~11.17]{GP}, whose Frobenius construction uses Larson--Sweedler \cite{LS}. Its generation and module-compatible depth identity are standard consequences of the integral and Hopf-module theorems \cite[Theorems~2.1.3 and~1.9.4]{Montgomery}. The adjunction comparison uses \cite[Lemmas~3.3--3.4]{Chen}; related comparisons up to retracts appear in \cite{AB}. The new reduction is the passage from one depth-two identity to the absorbing family for all objects. For unitary fusion categories, combining the characterization with \cite[Corollary~4.5]{EJ} gives an algebra-object criterion for on-site spin-chain symmetry; see Corollary~\ref{cor:onsite}.

\smallskip\noindent\emph{Organization.}
Section~\ref{sec:prelim} fixes the definitions, and Section~\ref{sec:depth} develops depth and its depth-two form. Section~\ref{sec:collapse} extracts a projective generator from the depth condition. Section~\ref{sec:fiber} constructs the fiber functor and identifies the distinguished algebra. Section~\ref{sec:hopf} applies the regular Hopf-algebra model and completes the proof. Section~\ref{sec:examples} records the Frobenius formulation, its application to on-site spin-chain symmetries, and a nonsemisimple example. 

\section{Definitions and elementary categorical facts}\label{sec:prelim}

\subsection{Finite categories and tensor categories}
A \emph{finite $k$-linear abelian category} is one equivalent, as a $k$-linear category, to the category of finite-dimensional left modules over a finite-dimensional associative unital $k$-algebra. In particular it has finite-dimensional morphism spaces, finite-length objects, finitely many simple isomorphism classes, and enough projective and injective objects. A nonzero object is \emph{simple} if its only subobjects are zero and itself. A \emph{composition series} of an object $X$ is a finite chain of subobjects $0=X_0\subset X_1\subset\cdots\subset X_n=X$ with each quotient $X_i/X_{i-1}$ simple. These quotients are its \emph{composition factors}; their isomorphism classes with multiplicities are independent of the series by the Jordan--H\"older theorem \cite[Theorem~1.5.4]{EGNO}.

An object $P$ is \emph{projective} when $\Hom(P,-)$ is exact, and an object $I$ is \emph{injective} when $\Hom(-,I)$ is exact. A \emph{projective generator} is a projective object $P$ such that every object is a quotient of a finite direct sum of copies of $P$. A \emph{projective cover} of a simple object $L$ is an epimorphism $P\to L$ from a projective object such that no proper subobject of $P$ maps onto $L$. Projective covers exist and are unique up to isomorphism in a finite category. A category is \emph{semisimple} if every object is a finite direct sum of simple objects; equivalently, every short exact sequence splits. We use the conventions of \cite[\S\S1.8 and 4.16]{EGNO}.

A \emph{monoidal category} has a tensor product, a unit, and associativity and unit isomorphisms satisfying the pentagon and triangle identities. It is \emph{rigid} if every object has left and right duals. We use a left dual $X^*$ with maps
\[
\ev_X:X^*\otimes X\to\one,
\qquad \coev_X:\one\to X\otimes X^*
\]
satisfying the two triangle identities. A right dual has evaluation in the opposite order. A finite tensor category is as specified in the introduction; exactness means exactness in each variable. We suppress associativity and unit isomorphisms using coherence.

A \emph{strong monoidal functor} $F$ is equipped with invertible natural maps
\[
F(X)\otimes F(Y)\longrightarrow F(X\otimes Y),
\qquad k\longrightarrow F(\one),
\]
compatible with associativity and units. A fiber functor must also be exact and faithful. In particular, the usual forgetful functor on representations of a quasi-Hopf algebra need not be a fiber functor: its tensor constraints may not exist.

For a finite-dimensional vector space $V$ and an object $M$ in a $k$-linear abelian category, the \emph{external tensor product} $V\boxt M$ is characterized by
\[
\Hom(V\boxt M,N)\cong\Hom_k(V,\Hom(M,N)).
\]
It is a direct sum of $\dim_k V$ copies of $M$, with no preferred choice of basis. External tensor products are used throughout to avoid confusing multiplicity spaces with objects of $\C$.

\subsection{Algebras, modules, and Frobenius structures}
An \emph{algebra object} $(A,m,e)$ consists of multiplication $m:A\otimes A\to A$ and unit $e:\one\to A$ satisfying associativity and the two unit identities. A \emph{right $A$-module} $(M,r_M)$ has a map $r_M:M\otimes A\to M$ satisfying
\[
r_M\circ(r_M\otimes\id_A)=r_M\circ(\id_M\otimes m),
\qquad r_M\circ(\id_M\otimes e)=\id_M.
\]
Module morphisms commute with these actions. The regular right module is $(A,m)$. The category $\C_A$ is abelian \cite[Exercise~7.8.7]{EGNO}; kernels and cokernels are computed in $\C$, with actions induced by exactness of tensoring by $A$. Its forgetful functor $U:\C_A\to\C$ is exact and faithful.

The \emph{free-module functor} is
\[
T:\C\longrightarrow\C_A,\qquad T(X)=X\otimes A.
\]
It is exact and left adjoint to $U$ \cite[Lemma~7.8.12]{EGNO}. Explicitly,
\begin{equation}\label{eq:freeadj}
\Hom_{\C_A}(X\otimes A,M)\cong\Hom_\C(X,U(M)).
\end{equation}
The forward map restricts along $\id_X\otimes e$; the inverse sends $f:X\to M$ to $r_M\circ(f\otimes\id_A)$. For $X=\one$ and $M=A$, this identifies $\End_{\C_A}(A)$ with $\Hom_\C(\one,A)$. Thus connectedness is exactly
\begin{equation}\label{eq:brick}
                       \End_{\C_A}(A)=k\id_A\cong k.
\end{equation}

\begin{samepage}
\begin{definition}\label{def:frobenius}
A \emph{Frobenius algebra object} is an algebra object together with a coassociative comultiplication $d:A\to A\otimes A$ and a counit $\lambda:A\to\one$ such that
\begin{equation}\label{eq:frobenius}
d\circ m=(m\otimes\id_A)\circ(\id_A\otimes d)
      =(\id_A\otimes m)\circ(d\otimes\id_A).
\end{equation}
A connected Frobenius algebra is \emph{separable} if $m\circ d=c\id_A$ for some $c\in k^\times$. Equivalently, $m$ has an $A$-bimodule section \cite[Definition~2.14 and equation~(2.15)]{FG}; an $A$-bimodule has commuting left and right actions. For an arbitrary algebra object, separability means this latter splitting condition. The equivalence with the scalar identity for the specified Frobenius structure requires connectedness.
\end{definition}
\end{samepage}

The Frobenius pairing is $b=\lambda\circ m:A\otimes A\to\one$. Its inverse copairing is $c=d\circ e:\one\to A\otimes A$; they satisfy
\begin{equation}\label{eq:snakes}
(b\otimes\id_A)\circ(\id_A\otimes c)=\id_A,
\qquad
(\id_A\otimes b)\circ(c\otimes\id_A)=\id_A.
\end{equation}
Thus $A$ is both a left and a right dual of itself. Conversely, an algebra with a functional $\lambda$ whose pairing $\lambda\circ m$ has an inverse copairing is Frobenius: set
\[
d=(m\otimes\id_A)\circ(\id_A\otimes c)
  =(\id_A\otimes m)\circ(c\otimes\id_A).
\]
This is the standard pairing characterization of Frobenius algebras; see \cite[Proposition~11.10]{GP} and \cite[\S3]{Muger}. Neither symmetry of $b$ nor a nonzero value of $\lambda\circ e$ is part of the definition.

\subsection{Tensor generation and additive closure}
A \emph{subquotient} of $Z$ means a quotient of a subobject of $Z$. The \emph{tensor-power generation} condition means that every object is a subquotient of some finite sum of powers of $A$, including $A^{\otimes0}=\one$. If $Z'\hookrightarrow Z$ and $Z'\twoheadrightarrow X$, exactness gives $Z'\otimes W\hookrightarrow Z\otimes W$ and $Z'\otimes W\twoheadrightarrow X\otimes W$. Thus tensoring on either side preserves subquotients; applying this twice proves that the tensor product of two subquotients is a subquotient of the tensor product. If $A$ is Frobenius, duality and \eqref{eq:snakes} show that the class of these subquotients is also closed under both duals; it is the full tensor subcategory generated by $A$ and closed under subquotients.

For an object $M$ in an additive category $\mathcal D$, its \emph{additive closure} $\add_{\mathcal D}(M)$ is the full subcategory of direct summands of finite direct sums of copies of $M$ \cite[\S2]{Krause}. More generally, $\add_{\mathcal D}\{M_i\}$ allows finite sums of objects from the specified family. A direct summand is also called a \emph{retract}: $N$ is a retract of $M$ when there are morphisms $i:N\to M$ and $p:M\to N$ with $p\circ i=\id_N$.

\section{Depth and the alternating tower}\label{sec:depth}

We first construct the two sectors and define depth for every algebra object. We then identify the depth-two condition and compare two-categorical realizations. Throughout, direct summands are taken in the stated Hom category.

\subsection{The two sectors and the canonical realization}
A \emph{bicategory} has objects, categories of $1$-morphisms between objects, and composition associative up to coherent isomorphism. Morphisms in the Hom categories are called $2$-morphisms. A $2$-category is a strict bicategory. We use additive $k$-linear Hom categories and composition additive in each variable. An additive category is \emph{idempotent complete} if every idempotent endomorphism has an image as a direct summand.

Suppose $Y:b\to a$ and $Y^\vee:a\to b$ are $1$-morphisms with an adjunction $Y^\vee\dashv Y$. Write its coevaluation and evaluation as
\[
\coev:\one_a\longrightarrow Y\circ Y^\vee,
\qquad \ev:Y^\vee\circ Y\longrightarrow\one_b.
\]
Thus coevaluation is the adjunction unit, and evaluation is its counit. They make
\[
A=Y\circ Y^\vee,\qquad
m=\id_Y\circ\ev\circ\id_{Y^\vee},\qquad e=\coev
\]
an algebra object in $\C=\operatorname{Hom}(a,a)$. By a \emph{realization of the algebra} $A$ we mean such an adjunction together with an identification of the induced multiplication and unit with those of $A$.

A reverse adjunction $Y\dashv Y^\vee$, with maps
$\widetilde{\coev}:\one_b\to Y^\vee\circ Y$ and
$\widetilde{\ev}:Y\circ Y^\vee\to\one_a$, supplies a Frobenius structure:
\[
d=\id_Y\circ\widetilde{\coev}\circ\id_{Y^\vee},
\qquad \lambda=\widetilde{\ev}.
\]
The triangle identities give the counit identities, while interchange gives coassociativity and the Frobenius identities \cite[Lemma~3.4]{Muger}. A single adjunction supplies only the algebra structure. We require no adjoints for the other $1$-morphisms of the bicategory. M\"uger's Frobenius construction addresses two-sided duality and, with additional hypotheses, Morita equivalence; the one-sided realization used here addresses the algebra structure alone. See also \cite[Proposition~2.26]{GP}.

Put $\M=\operatorname{Hom}(b,a)$. We call $\C$ the \emph{even sector} and $\M$ the \emph{odd sector}; a sector here is a Hom category. The associated levels are
\begin{equation}\label{eq:abstract-levels}
 E_r=A^{\otimes r}\in\C\quad\text{at level }2r,
 \qquad A^{\otimes r}\circ Y\in\M\quad\text{at level }2r+1.
\end{equation}
An odd summand is compared only with earlier odd levels, and an even summand only with earlier even levels.

Every algebra in a finite tensor category has a canonical realization in the bicategory of algebras and bimodules internal to $\C$. Restrict to its objects $\one$ and $A$. Composition through $A$ is the relative tensor product, defined as the coequalizer of the two action maps. Unlike the ambient tensor product, a relative tensor product is generally only right exact \cite[Exercise~7.8.23]{EGNO}. The even category at $\one$ is $\C$, and the adjacent category is $\C_A$. Take
\[
Y_A={}_{\one}A_A,\qquad Y_A^\vee={}_AA_{\one}.
\]
Their even composite is $A\otimes_A A\cong A$. The coevaluation of $Y_A^\vee\dashv Y_A$ is $e:\one\to A$, and its evaluation is $m:A\otimes A\to A$, viewed as an $A$-bimodule map. After the relative-tensor unit identifications, its triangle composites are
\[
A\xrightarrow{e\otimes\id_A}A\otimes A\xrightarrow{m}A,
\qquad
A\xrightarrow{\id_A\otimes e}A\otimes A\xrightarrow{m}A;
\]
both are the identity by the unit identities of $A$. A Frobenius structure supplies the reverse adjunction with evaluation $\lambda:A\to\one$ and coevaluation $d:A\to A\otimes A$ in the $A$-bimodule sector. Conversely, a reverse adjunction on this canonical pair supplies a Frobenius structure. Thus no Frobenius hypothesis is needed for the first adjunction or for the definition below. For the bimodule realization, see \cite[Remark~3.18]{Muger}.

The \emph{canonical sectors} are consequently $\C$ and $\C_A$, and their levels are
\begin{equation}\label{eq:canonical-levels}
 E_r=A^{\otimes r}\in\C,
 \qquad O_r=A^{\otimes r}\otimes A\in\C_A,
\end{equation}
where only the last factor carries the right action. In particular, level one is the regular right module $A$, not the unit object of $\C$.

\subsection{Definition of depth}
An object is \emph{indecomposable} if it is nonzero and is not a direct sum of two nonzero objects. The category $\C_A$ is finite abelian by \cite[Exercises~7.8.7 and~7.8.16]{EGNO}. Hence both $\C$ and $\C_A$ have the Krull--Schmidt property: every object is a finite direct sum of indecomposables, uniquely up to isomorphism and permutation \cite[Theorem~1.5.7]{EGNO}.

Define the cumulative additive closures
\[
 \mathcal E_r=\add_{\C}\!\left(\bigoplus_{j=0}^r E_j\right),
 \qquad
 \mathcal O_r=\add_{\C_A}\!\left(\bigoplus_{j=0}^r O_j\right),
 \qquad \mathcal E_{-1}=\mathcal O_{-1}=\{0\}.
\]
\begin{samepage}
\begin{definition}[Depth]\label{def:all-depth}
The \emph{depth} of a nonzero algebra object $A$ is
\[
\depth(A)=\sup\Bigl(
 \{2r:\mathcal E_r\ne\mathcal E_{r-1}\}
 \ \cup\
 \{2r+1:\mathcal O_r\ne\mathcal O_{r-1}\}
\Bigr).
\]
Equivalently, depth is the last level at which a new isomorphism class of indecomposable direct summand first appears in its own sector. It is infinity if these first appearances are unbounded. The algebra has \emph{finite depth} precisely when this supremum is finite.
\end{definition}
\end{samepage}

\begin{example}[A finite group]\label{ex:finite-group}
Let $G$ be a finite group, let $\C=\Vect_G$ be the category of
finite-dimensional $G$-graded $k$-vector spaces, and let $A=kG$
with its usual grading.
For every graded right $A$-module $M$, multiplication by $g$
identifies $M_e$ with $M_g$. Thus $M\mapsto M_e$ gives an equivalence
$\C_A\simeq\Vect_k$, under which the regular module $A_A$
corresponds to $k$.

Write $k_g$ for the line in degree $g$. Level~$0$ contains only
$k_e$. At level~$1$, $A_A$ is the unique odd indecomposable.
At level~$2$, $A=\bigoplus_{g\in G}k_g$ contains every even
indecomposable. Thus, for $G\ne\{e\}$, the last new summands appear
at level~$2$, so $\depth(A)=2$ by Definition~\ref{def:all-depth}.
For the trivial group, the depth is~$1$.
\end{example}

The cumulative formulation does not require the unit to split off $A$. Finiteness of depth is equivalent to finiteness of the set of summand classes appearing in the two towers: a finite set has a largest first-appearance level, and finitely many levels contain only finitely many summand classes. A finite tensor category can have infinitely many indecomposable objects, so its finiteness does not by itself imply finite depth. For example, over an infinite field of characteristic two, the $C_2\times C_2$-modules on $k^2$ whose generators act as $I+N$ and $I+aN$, where $N=\left(\begin{smallmatrix}0&1\\0&0\end{smallmatrix}\right)$ and $a\in k$, are pairwise nonisomorphic indecomposables: the first generator has one Jordan block, and an intertwining isomorphism forces the parameters to agree.

When both sectors are semisimple, indecomposables are simple and simple subquotients are direct summands. Thus $\mathcal E_r$ records exactly the simple classes appearing up to level $2r$, and $\mathcal O_r$ those appearing up to level $2r+1$ in the odd sector. Their change levels are precisely the first-appearance levels in the alternating tower rooted at $\one$. Proposition~\ref{prop:subfactor-depth} identifies this convention with principal-graph depth for subfactors.

\begin{samepage}
\begin{remark}\label{rem:composition-factor-depth}
Replacing indecomposable direct summands by simple subquotients in the alternating tower defines \emph{composition-factor depth}. The two depths agree when $\C$ and $\C_A$ are semisimple, in particular when $\C$ is semisimple and $A$ is separable \cite[Proposition~7.8.30]{EGNO}. The distinction is essential: for the complex quasi-Hopf algebra $H_+(8)$ of \cite{EG8}, $\Rep(H_+(8))$ has no fiber functor, yet is tensor-generated by a connected separable Frobenius algebra of composition-factor depth two and depth five.
\end{remark}
\end{samepage}

\subsection{The depth-two identity}
\begin{samepage}
\begin{lemma}[Summands of copies of an object with scalar endomorphisms]\label{lem:matrix}
Let $M$ be an object of an additive $k$-linear category with $\End(M)=k$. Every object of $\add(M)$ is isomorphic to $M^{\oplus s}$ for some $s\geq0$.
\end{lemma}
\end{samepage}
\begin{proof}
By \cite[Proposition~2.3]{Krause}, $\Hom(M,-)$ is fully faithful on $\add(M)$ and takes values in finitely generated projective right $\End(M)$-modules. Here $\End(M)=k$, so every such module is $k^{\oplus s}=\Hom(M,M^{\oplus s})$. Full faithfulness gives the required isomorphism.
\end{proof}

\begin{samepage}
\begin{proposition}\label{prop:depth-two-identity}
For a nonzero algebra object $A$,
\begin{equation}\label{eq:depth-module}
\depth(A)\leq2
\quad\Longleftrightarrow\quad
A\otimes A\in\add_{\C_A}(A),
\end{equation}
where the right action on $A\otimes A$ is $\id_A\otimes m$. Equivalently, in any idempotent-complete realization as in \eqref{eq:abstract-levels},
\[
                      A\circ Y\in\add_{\M}(Y).
\]
If $A$ is connected, these conditions are equivalent to
\[
 A\otimes A\cong A^{\oplus n}\ \text{in }\C_A,
 \qquad\text{or, respectively,}\qquad
 A\circ Y\cong Y^{\oplus n}
\]
for some integer $n\geq1$.
\end{proposition}
\end{samepage}
\begin{proof}
If the depth is at most two, level three adds no indecomposable summand, so $O_1\in\add(O_0)$. Conversely, suppose this condition holds. The functor $A\otimes-:\C_A\to\C_A$, with action $\id_A\otimes r_M$ on $A\otimes M$, is additive. It preserves finite sums by the biproduct identities. It preserves a retraction $p\circ i=\id_M$ because
\[
(\id_A\otimes p)\circ(\id_A\otimes i)=\id_{A\otimes M}.
\] Induction gives $O_r\in\add(O_0)$ for every $r\geq0$. Forgetting the module structure gives $E_r\in\add_{\C}(A)$ for all $r\geq1$. Thus the odd additive closures are constant from level one and the even closures from level two. There are no later new summands.

Proposition~\ref{prop:comparison} identifies the additive closures in any realization with the canonical ones, giving the equivalent condition on $Y$. For connected $A$, \eqref{eq:brick} gives $\End_{\C_A}(A)=k$, and the comparison also gives $\End_{\M}(Y)=k$. Lemma~\ref{lem:matrix} therefore applies to both $A$ and $Y$. Finally $A\otimes A\ne0$: its subobject $\one\otimes A$, induced by the monomorphism $e$, is nonzero. Hence $n\geq1$.
\end{proof}

In concrete terms, \eqref{eq:depth-module} asks for right-module morphisms
\[
A\otimes A\xto{i} A^{\oplus r}\xto{p} A\otimes A,
\qquad p\circ i=\id_{A\otimes A}.
\]
The requirement that these maps commute with right multiplication is part of the condition. Left multiplication gives the analogous left-module condition. Theorem~\ref{thm:main} will imply that the left condition holds as well for every connected tensor-generating algebra of depth at most two.

\subsection{Comparison with two-categorical realizations}
Definition~\ref{def:all-depth} is intrinsic: it uses $A$ in $\C$ and its regular right module in $\C_A$. No choice of factorization is required. The purpose of the following comparison is to show that any adjunction realizing the \emph{algebra} $A$ computes this same invariant from its own alternating tower. In particular it justifies the formulation using $Y$ in Proposition~\ref{prop:depth-two-identity}. The argument needs only the one adjunction, without Frobenius or separability assumptions; its mechanism is standard \cite[Lemmas~3.3--3.4]{Chen}.

\begin{samepage}
\begin{proposition}[Comparison on the additive tower]\label{prop:comparison}
Consider a factorization as in \eqref{eq:abstract-levels}, with idempotent-complete Hom categories. The canonical comparison functor induces an equivalence
\[
 \add_{\M}\{Z\circ Y:Z\in\C\}
 \ \simeq\
 \add_{\C_A}\{Z\otimes A:Z\in\C\},
\]
sending $A^{\otimes r}\circ Y$ to $O_r$. It induces the corresponding equivalence on the additive closures of each finite list of odd levels. Hence depth is independent of the chosen realization of the algebra.
\end{proposition}
\end{samepage}
\begin{proof}
Put $F(Z)=Z\circ Y$ and $R(M)=M\circ Y^\vee$. The given adjunction induces $F\dashv R$, with monad $R\circ F=-\otimes A$ and the prescribed algebra structure. Its comparison functor $K:\M\to\C_A$ sends $F(Z)$ to $Z\otimes A$ and is fully faithful on these objects by \cite[Lemma~3.3]{Chen}.

Let $\mathcal P$ and $\mathcal Q$ be the two additive closures in the statement, and let $\mathcal P_0$ be the full subcategory of finite sums of objects $F(Z)$. In the notation of \cite[Lemma~3.4]{Chen}, take
\[
(F,\mathcal C,\mathcal C',\mathcal D)
=(K|_{\mathcal P},\mathcal P,\mathcal P_0,\mathcal Q).
\]
The object denoted $Y$ there is an arbitrary module $N\in\mathcal Q$, which is a retract of $K(\bigoplus_iF(Z_i))$. Likewise, every object of $\mathcal P$ is a retract of an object of $\mathcal P_0$. Additivity gives full faithfulness on $\mathcal P_0$, so parts~(1)--(2) of that lemma give an equivalence up to retracts. Since $\mathcal P$ is idempotent complete, part~(3) gives an equivalence.

For a finite list of odd levels, restrict $\mathcal P_0$ to sums of $F(A^{\otimes r})$ from that list and take their retract closures on both sides. The same argument preserves each level and its indecomposable summands, hence their first appearances. The even levels are unchanged.
\end{proof}

The equivalence identifies the additive closures of the entire tower, as well as those at each finite stage. It therefore preserves both finite depth and the value of depth, including infinity.

\subsection{Agreement with subfactor depth}\label{subsec:subfactor-depth}
We use the principal graph rooted at the trivial $N$--$N$ bimodule, with vertex depth given by distance from that root; see \cite[\S2.1.3]{JMS} and \cite[Definition~2.1]{MS}.

\begin{samepage}
\begin{proposition}[Agreement with principal-graph depth]\label{prop:subfactor-depth}
Let $N\subset M$ be a finite-index, finite-depth inclusion of $\mathrm{II}_1$ factors. Let $A={}_NL^2(M)_N$ be its canonical algebra in the generated $N$--$N$ fusion category, and let $\Gamma_N$ be its principal graph rooted at $\one_N=L^2(N)$. Then
\begin{equation}\label{eq:subfactor-depth}
\depth(A)=\max_{v\in\Gamma_N}\operatorname{dist}(\one_N,v).
\end{equation}
\end{proposition}
\end{samepage}
\begin{proof}
Put $Y={}_NL^2(M)_M$ and $Y^\vee={}_ML^2(M)_N$, with their standard conjugate adjunction. Its induced algebra is $A\cong Y\otimes_M Y^\vee$. The alternating levels are $W_{2r}=A^{\otimes_N r}$ and $W_{2r+1}=A^{\otimes_N r}\otimes_NY$. Proposition~\ref{prop:comparison} identifies their odd additive closures with those of the free right modules $O_r$; the even levels coincide. By the standard principal-graph description \cite[\S2.1.3]{JMS}, edges of $\Gamma_N$ record tensoring with $Y$, and adjunction gives the reverse edge multiplicities by tensoring with $Y^\vee$. Inductively, the multiplicity of a vertex $v$ in $W_n$ is the number of length-$n$ edge walks from the root to $v$. Its first appearance is therefore at its graph distance. Semisimplicity identifies these first appearances with the changes of the additive closures in Definition~\ref{def:all-depth}.
\end{proof}

For an infinite-depth finite-index inclusion, the generated even category is not a finite tensor category; extending the same tower formula to its semisimple finite-length sectors gives the same equality with supremum and value infinity. The dual algebra computes the $M$-rooted depth $d_M$, which can differ from $d_N=\depth(A)$ by one \cite[\S2.1]{MS}; the maximum-of-the-pair convention gives $\max\{d_N,d_M\}$. The identity inclusion has depth one.

\section{From depth two to a projective generator}\label{sec:collapse}

The main step is a statement about all modules, obtained without assuming that $\C_A$ is semisimple. By Proposition~\ref{prop:depth-two-identity}, the depth hypothesis is exactly \eqref{eq:depth-module}. We use the fact that projective and injective objects coincide in a finite tensor category \cite[Proposition~6.1.3 and Remark~6.1.4]{EGNO}.

\begin{samepage}
\begin{proposition}[The module category collapses]\label{prop:collapse}
Let $A$ satisfy condition \textup{(ii)} of Theorem~\ref{thm:main}. Then $A$ is a projective generator of $\C_A$, and every object of $\C_A$ is a finite direct sum of copies of $A$. The evaluation maps
\begin{equation}\label{eq:evaluation}
e_M:\Hom_{\C_A}(A,M)\boxt A\longrightarrow M
\end{equation}
are natural isomorphisms. Consequently $V\mapsto V\boxt A$ is an equivalence $\Vect_k\simeq\C_A$.
\end{proposition}
\end{samepage}

\begin{proof}
We divide the argument into four steps. Evaluation in \eqref{eq:evaluation} is defined by the external-tensor universal property, as detailed in Step~4.

\smallskip\noindent\emph{Step 1: tensor powers reduce to copies of $A$.}
Connectedness gives \eqref{eq:brick}. Lemma~\ref{lem:matrix} turns the depth condition into
\begin{equation}\label{eq:module-power}
A\otimes A\cong A^{\oplus n}\quad\text{as right $A$-modules}
\end{equation}
for some $n>0$. Indeed, $A\neq0$, and $A\otimes A\neq0$ since $e:\one\to A$ is a monomorphism and $A\otimes-$ is exact. Tensoring \eqref{eq:module-power} on the left and proceeding by induction gives
\begin{equation}\label{eq:all-powers}
A^{\otimes r}\cong A^{\oplus n^{r-1}}
\quad\text{as right $A$-modules},\qquad r\geq1,
\end{equation}
with the action always on the last tensor factor. Moreover $\one$ is a subobject of $A$ via $e$. Generation now implies that every object of $\C$ is a subquotient of a finite direct sum of copies of the underlying object $A$.

\smallskip\noindent\emph{Step 2: a projective generator is a summand of copies of $A$.}
Choose a projective generator $P$ of $\C$. By Step~1 there is a subobject $Z\subseteq A^{\oplus r}$ and an epimorphism $Z\to P$. Projectivity splits this epimorphism, so $P$ embeds in $A^{\oplus r}$. By \cite[Proposition~6.1.3 and Remark~6.1.4]{EGNO}, $P$ is also injective. Its embedding therefore splits, proving
\begin{equation}\label{eq:projective-summand}
                         P\in\add_\C(A).
\end{equation}
This is the step that recovers information about extensions from the tensor-generation hypothesis.

\smallskip\noindent\emph{Step 3: free modules produce a projective generator in $\add(A)$.}
The adjunction \eqref{eq:freeadj}, together with exactness of $U$, shows that
$Q=T(P)=P\otimes A$ is projective in $\C_A$. It is also a generator. To see this, choose an epimorphism $P^{\oplus t}\to U(M)$. Tensor it with $A$ and then apply the action:
\[
(P\otimes A)^{\oplus t}\longrightarrow U(M)\otimes A
\xto{r_M}M.
\]
Both maps are epimorphisms; the second has a section in $\C$ given by $\id_M\otimes e$. Thus every module is a quotient of copies of $Q$.

Applying $T$ to the splitting in \eqref{eq:projective-summand} yields
\[
Q\in\add_{\C_A}(A\otimes A)\subseteq\add_{\C_A}(A).
\]
Lemma~\ref{lem:matrix} gives $Q\cong A^{\oplus s}$ for some $s>0$. It follows that $A$ itself is a projective generator of $\C_A$.

\smallskip\noindent\emph{Step 4: apply the projective-generator equivalence.}
By the projective-generator form of Morita theory \cite[\S1.8, pp.~9--10]{EGNO},
\[
G=\Hom_{\C_A}(A,-):\C_A\xrightarrow{\sim}R\text{-}\mathrm{mod},
\qquad R=\End_{\C_A}(A)^{\op}=k.
\]
Its inverse is $L(V)=V\boxt A$, and its counit is the evaluation map $e_M$ in \eqref{eq:evaluation}. In particular every $M$ is a sum of copies of $A$ and every $e_M$ is invertible. For later use, we record the standard evaluation identities. If $\iota_f:A\to G(M)\boxt A$ is induced by $1\mapsto f$, then
\begin{equation}\label{eq:eval-on-copower}
e_M\circ\iota_f=f.
\end{equation}
With $G(h)(f)=h\circ f$, counit naturality says
\begin{equation}\label{eq:eval-naturality}
h\circ e_M=e_N\circ(G(h)\boxt\id_A).
\end{equation}
For any isomorphism $\alpha:M\xrightarrow{\sim}A^{\oplus d}$ this gives
\begin{equation}\label{eq:eval-transport}
e_M=\alpha^{-1}\circ e_{A^{\oplus d}}\circ(G(\alpha)\boxt\id_A).
\end{equation}
\end{proof}

\begin{samepage}
\begin{remark}[From one identity to the absorbing family]\label{rem:absorption}
For every $X\in\C$, Proposition~\ref{prop:collapse} applied to $M=X\otimes A$ gives the natural isomorphism
\[
\Hom_{\C_A}(A,X\otimes A)\boxt A
 \xrightarrow{\ e_{X\otimes A}\ }X\otimes A.
\]
Thus the single depth-two identity supplies the absorbing family for all objects. The passage beyond tensor powers occurs in Steps~2--3: tensor generation puts a projective generator of $\C$ in $\add_\C(A)$, and the free-module functor puts a projective generator of $\C_A$ in $\add_{\C_A}(A)$. Step~4 gives the displayed decomposition. Simplicity of $A$ as a module and semisimplicity of $\C_A$ are conclusions of this argument.
\end{remark}
\end{samepage}

\section{The fiber functor and the distinguished algebra}\label{sec:fiber}

\subsection{Tensor structure and exactness}
Assume that $L:\Vect_k\to\C_A$, $L(V)=V\boxt A$, is an equivalence, as Proposition~\ref{prop:collapse} ensures under the depth-two hypotheses. Its quasi-inverse is
\[
G(M)=\Hom_{\C_A}(A,M),\qquad G(h)(f)=h\circ f.
\]
The unit $j_V:V\to G(L(V))$, $v\mapsto\iota_v$, and the evaluation counit $e_M:L(G(M))\to M$ of the external-tensor adjunction are isomorphisms. Thus $G\circ L\cong\id_{\Vect_k}$ and $L\circ G\cong\id_{\C_A}$. Define
\begin{equation}\label{eq:fiber}
F_A(X)=\Hom_{\C_A}(A,X\otimes A)
       \cong\Hom_\C(\one,X\otimes A).
\end{equation}
On a morphism $u:X\to Y$, the map $F_A(u)$ is postcomposition with $u\otimes\id_A$.

\begin{samepage}
\begin{proposition}\label{prop:fiber}
The tensor constraints
\begin{align}
\mu_{X,Y}:F_A(X)\otimes_k F_A(Y)&\longrightarrow F_A(X\otimes Y),\label{eq:mu}\\[-6pt]
f\otimes g&\longmapsto(\id_X\otimes g)\circ f\notag
\end{align}
and $k\to F_A(\one)$, $1\mapsto\id_A$, make $F_A$ an exact faithful strong monoidal functor.
\end{proposition}
\end{samepage}

\begin{proof}
Transport the exact left $\C$-action on $\C_A$ through $L$. It becomes the action
\[
X\mathbin{\triangleright}V=G(X\otimes L(V))
\]
on $\Vect_k$. The standard correspondence between such actions and fiber functors \cite[Proposition~7.3.3 and Example~7.4.6]{EGNO} gives the functor $X\mapsto X\mathbin{\triangleright}k=F_A(X)$. Under the evaluation identifications, the action associativity sends $f\otimes g$ to $(\id_X\otimes g)\circ f$, and its unit sends $1$ to $\id_A$. These are exactly the stated constraints. Exactness follows from exactness of $X\mapsto X\otimes A$ and of the inverse equivalence $G$; faithfulness also follows from \cite[Remark~4.3.10]{EGNO}.
\end{proof}

\subsection{Why reconstruction carries the algebra to the regular comodule}
Finite-dimensional reconstruction gives a Hopf algebra $B=\Nat(F_A,F_A)$ and a tensor equivalence $\C\simeq\Rep(B)$ whose underlying vector-space functor is $F_A$; see \cite[Theorem~5.3.12]{EGNO} and \cite[\S\S2.2--2.3]{Schauenburg} for reconstruction over a field. With our left-comodule convention, put $H=(B^*)^{\mathrm{cop}}$ and denote the resulting tensor equivalence by
\[
E:\C\xrightarrow{\sim}\Corep(H).
\]
The coopposite leaves multiplication and unit unchanged and has
\[
\Delta_H(h)(b\otimes c)=h(cb),\qquad S_H(h)=h\circ S_B^{-1}.
\]
The inverse antipode exists in finite dimension \cite[Proposition~5.3.15]{EGNO}.

\begin{samepage}
\begin{proposition}\label{prop:regular}
Under this equivalence, $A$ becomes $H$ as a left comodule \emph{algebra}.
\end{proposition}
\end{samepage}
\begin{proof}
Write $\underline{\End}$ for the internal endomorphism algebra of an object in a module category. The free-module adjunction identifies $\underline{\End}_{\C_A}(A)$ with $(A,m,e)$: evaluation is $m$, and internal composition and its unit are the given multiplication and unit. The module equivalence $L$, with $L(k)=A$, therefore gives
\[
E(A)\cong\underline{\End}_{\Vect_k}(k)\cong H
\]
as algebra objects. The last identification is \cite[Examples~7.9.10--7.9.11]{EGNO}: in $\Rep(B)$ it is $B^*$ with right translations $(b\cdot h)(c)=h(cb)$, which is the regular left comodule algebra for $H=(B^*)^{\mathrm{cop}}$.
\end{proof}

The Frobenius structure and the remaining structural conclusions follow from the regular Hopf-algebra model, which we treat next.

\section{The regular Hopf-algebra model}\label{sec:hopf}

\subsection{Hopf algebras and left comodules}
A \emph{bialgebra} over $k$ is an associative unital algebra $H$ with a coassociative coproduct $\Delta_H:H\to H\otimes_k H$ and counit $\eps_H:H\to k$, both algebra homomorphisms. A \emph{Hopf algebra} has an antipode $S:H\to H$ satisfying
\[
m_H\circ(S\otimes\id)\circ\Delta_H
=e_H\circ\eps_H
=m_H\circ(\id\otimes S)\circ\Delta_H.
\]
For finite-dimensional Hopf algebras the antipode is bijective \cite[Proposition~5.3.15]{EGNO}. We use Sweedler notation $\Delta_H(h)=\sum h_{(1)}\otimes h_{(2)}$, suppressing summation signs in computations.

A \emph{left $H$-comodule} is a vector space $X$ with a map
$\rho_X:X\to H\otimes X$, written $x\mapsto x_{(-1)}\otimes x_{(0)}$, such that
\[
(\Delta_H\otimes\id_X)\circ\rho_X=(\id_H\otimes\rho_X)\circ\rho_X,
\qquad(\eps_H\otimes\id_X)\circ\rho_X=\id_X.
\]
Its tensor product has coaction
\[
\rho_{X\otimes Y}(x\otimes y)
=x_{(-1)}y_{(-1)}\otimes x_{(0)}\otimes y_{(0)}.
\]
Finite-dimensional left comodules form $\Corep(H)$, a finite tensor category; equivalently they are finite-dimensional right $H^*$-modules. The regular left comodule is $(H,\Delta_H)$. Its multiplication and unit are colinear, so it is an algebra object. These are the conventions of \cite[Definition~11.4]{GP}; its Theorem~11.11 and Lemma~11.17 therefore apply directly over $k$.

\begin{samepage}
\begin{proposition}\label{prop:hopf-depth}
For every finite-dimensional Hopf algebra $H$ over $k$, its regular comodule algebra is connected, its tensor powers generate $\Corep(H)$, and
\[
H\otimes H\cong H^{\oplus\dim_k H}
\quad\text{as right $H$-modules in }\Corep(H).
\]
\end{proposition}
\end{samepage}

\begin{proof}
Connectedness follows from \cite[Theorem~11.11 and Lemma~11.17]{GP}. By the Larson--Sweedler theorem, $H$ is cyclic as a right $H^*$-module \cite[Theorem~2.1.3(3)]{Montgomery}. Since $\dim_kH=\dim_kH^*$, it is isomorphic to the regular right $H^*$-module. Under $\Corep(H)\simeq\operatorname{mod}\!\text{-}H^*$, this makes $H$ a projective generator; in particular it tensor-generates.

Right $H$-modules in $\Corep(H)$ are left-comodule/right-module Hopf modules. The fundamental theorem of Hopf modules, in this mixed-handed form, makes every such object a finite copower of $H$ \cite[Theorem~1.9.4 and the following paragraph, pp.~15--16]{Montgomery}; bijectivity of the antipode ensures that this form applies. Apply it to $H\otimes H$, with its tensor-product coaction and right multiplication on the second factor. Comparing vector-space dimensions gives the multiplicity $\dim_kH$.
\end{proof}

\subsection{The Frobenius structure and its two counits}
The Frobenius assertion is the following result of Ghosh--Palcoux.

\begin{samepage}
\begin{proposition}[Ghosh--Palcoux]\label{prop:hopf-frob}
The regular left comodule algebra $H$ is connected and Frobenius. The space of functionals $\lambda:H\to k$ satisfying
\begin{equation}\label{eq:invariant-integral}
h_{(1)}\lambda(h_{(2)})=\lambda(h)1
\end{equation}
is one-dimensional, and every nonzero such functional is a Frobenius counit for this algebra object.
\end{proposition}
\end{samepage}

See \cite[Theorem~11.11 and Lemma~11.17]{GP}; the passage from the nondegenerate integral pairing to a Frobenius structure is \cite[Proposition~11.10]{GP}. These results require neither semisimplicity nor algebraic closure.

For use in Example~\ref{ex:charp}, we record the pairing description from Section~\ref{sec:prelim}. The pairing
\begin{equation}\label{eq:hopf-pairing}
b_H(x,y)=\lambda(xy)
\end{equation}
has an inverse copairing $c=\sum_i u_i\otimes v_i$, characterized by
\begin{equation}\label{eq:dual-bases}
\sum_i u_i\lambda(v_i h)=h,
\qquad \sum_i\lambda(h u_i)v_i=h.
\end{equation}
The Frobenius coproduct is
\begin{equation}\label{eq:frob-coprod}
d(h)=\sum_i h u_i\otimes v_i=\sum_i u_i\otimes v_i h.
\end{equation}

\begin{samepage}
\begin{remark}\label{rem:two-coproducts}
The Hopf coproduct $\Delta_H$ defines the comodule structure, whereas $d$ is a morphism into the tensor-product comodule with the Frobenius identities. Likewise, the Hopf counit $\eps_H$ and the Frobenius counit $\lambda$ have different roles. In particular, one must not assume that $\eps_H:H\to k$ is a comodule morphism, or that $\lambda(1)$ is nonzero. The example in Section~\ref{sec:examples} makes this distinction explicit.
\end{remark}
\end{samepage}

\subsection{Proof of the characterization and structural consequences}
\begin{proof}[Proof of Theorem~\ref{thm:main}]
Given a fiber functor, reconstruction \cite[Theorem~5.3.12]{EGNO} identifies $\C$ with $\Corep(H)$ for a finite-dimensional Hopf algebra $H$. Proposition~\ref{prop:hopf-depth}, together with Proposition~\ref{prop:depth-two-identity}, supplies a connected tensor-generating algebra of depth at most two. This proves \textup{(i)}$\Rightarrow$\textup{(ii)}.

Conversely, let $A$ satisfy \textup{(ii)}. Proposition~\ref{prop:depth-two-identity} gives the right-module identity \eqref{eq:depth-module}. Proposition~\ref{prop:collapse} shows that every right $A$-module is a finite direct sum of copies of $A$. Proposition~\ref{prop:fiber} then constructs a fiber functor $F_A$, proving \textup{(ii)}$\Rightarrow$\textup{(i)}. Under reconstruction for this functor, Proposition~\ref{prop:regular} identifies $A$ with the regular comodule algebra. Proposition~\ref{prop:hopf-frob} supplies its Frobenius structure.
\end{proof}

\begin{samepage}
\begin{corollary}[The distinguished algebra]\label{cor:structure}
Let $A$ be a connected algebra that tensor-generates $\C$ and has depth at most two. Then:
\begin{enumerate}
\item $V\mapsto V\boxt A$ is an equivalence $\Vect_k\simeq\C_A$, with inverse $M\mapsto\Hom_{\C_A}(A,M)$.
\item Formula~\eqref{eq:fiber}, with the tensor constraints \eqref{eq:mu}, defines a fiber functor. Its reconstructed tensor equivalence $\C\simeq\Corep(H)$ sends $A$ to the regular comodule algebra $H$.
\item $A$ is a projective generator of $\C$ and admits a Frobenius structure. Moreover,
\begin{equation}\label{eq:dimension}
 A\otimes A\cong A^{\oplus n}\quad\text{in }\C_A,
 \qquad n=\dim_k H=\dim_k\End_{\C}(A).
\end{equation}
The corresponding left-module identity also holds.
\end{enumerate}
Conversely, for any algebra $A$ such that $V\mapsto V\boxt A$ is an equivalence $\Vect_k\simeq\C_A$, the algebra $A$ is connected, tensor-generates $\C$, and has depth at most two.
\end{corollary}
\end{samepage}
\begin{proof}
The equivalence and the construction of the fiber functor were proved in Propositions~\ref{prop:collapse} and \ref{prop:fiber}. Proposition~\ref{prop:regular} identifies $A$ with $H$, and Proposition~\ref{prop:hopf-frob} gives the Frobenius structure.

The proof of Proposition~\ref{prop:hopf-depth} identifies $H$, under $\Corep(H)\simeq\operatorname{mod}\!\text{-}H^*$, with the regular right $H^*$-module. Therefore $H$ is a projective generator and $\dim_k\End_{\Corep(H)}(H)=\dim_kH$. Comparing vector-space dimensions in $H\otimes H\cong H^{\oplus n}$ gives $n=\dim_kH$ as an equality of positive integers. The left-module identity follows by applying the left-comodule/left-module form of the fundamental theorem of Hopf modules to $H\otimes H$, with left multiplication on the first factor \cite[Theorem~1.9.4 and the following paragraph]{Montgomery}.

Finally, assume only the equivalence $\Vect_k\simeq\C_A$ in the converse. Proposition~\ref{prop:fiber} constructs $F_A$ directly from it, and Proposition~\ref{prop:regular} identifies $A$ with the regular comodule algebra of the reconstructed Hopf algebra. Proposition~\ref{prop:hopf-depth} then gives connectedness, tensor generation, and depth at most two.
\end{proof}

\begin{samepage}
\begin{corollary}[Automatic Frobenius structure]\label{cor:automatic-frobenius}
Every connected algebra object that tensor-generates a finite tensor category and has depth at most two admits a Frobenius structure.
\end{corollary}
\end{samepage}
\begin{proof}
Corollary~\ref{cor:structure} identifies it as an algebra with the regular left comodule algebra of a finite-dimensional Hopf algebra. Proposition~\ref{prop:hopf-frob} supplies the Frobenius structure.
\end{proof}

\section{Frobenius algebras of depth two}\label{sec:examples}

\subsection{The Frobenius formulation and the trivial case}
\begin{samepage}
\begin{corollary}\label{cor:frobenius}
A finite tensor category admits a fiber functor if and only if it is tensor-generated by a connected Frobenius algebra of depth at most two.
\end{corollary}
\end{samepage}
\begin{proof}
Theorem~\ref{thm:main} applies after forgetting the Frobenius structure. In the other direction, the same theorem supplies the Frobenius structure for an algebra satisfying its condition \textup{(ii)}.
\end{proof}

\begin{samepage}
\begin{corollary}\label{cor:exact-depth}
A connected tensor-generating algebra has depth at most one if and only if it is the unit algebra and $\C\simeq\Vect_k$. Consequently, in a nontrivial finite tensor category, every algebra in Theorem~\ref{thm:main}\textup{(ii)} has depth exactly two.
\end{corollary}
\end{samepage}
\begin{proof}
Depth at most one means that even level two introduces no new summand, so $A\in\add_{\C}(\one)$. Since $\End(\one)=k$, Lemma~\ref{lem:matrix} gives $A\cong\one^{\oplus r}$. Connectedness forces $r=1$, and the unit then identifies $A$ with the unit algebra. Tensor generation implies that every object is a subquotient of a sum of copies of the simple object $\one$. Such sums and their subquotients are semisimple, so all objects are sums of $\one$ and $\C\simeq\Vect_k$ as a tensor category. Conversely, for the unit algebra in $\Vect_k$, the only first appearances are $\one$ at level zero and the single odd class at level one. Its depth is one.
\end{proof}

\subsection{On-site symmetries of spin chains}\label{sec:onsite}
Let $\C$ be a unitary fusion category over $\mathbb C$, that is, a finite semisimple rigid $C^*$-tensor category with simple unit and unitary associativity and unit constraints. We use the categorical symmetries of Evans--Jones \cite{EJ}. The tensor-product quasi-local algebra
\[
\mathfrak A_q=\bigotimes_{\mathbb Z}M_q(\mathbb C)
\]
is the norm completion of the union of its finite-interval matrix algebras. A categorical symmetry is \emph{on-site} when each symmetry channel preserves every such interval algebra \cite[Definition~3.25]{EJ}. A \emph{Q-system} is a Frobenius algebra compatible with the adjoint, with $d=m^*$, $\lambda=e^*$, and $m\circ m^*$ a positive scalar multiple of $\id_A$.

\begin{samepage}
\begin{corollary}[On-site realizability]\label{cor:onsite}
The category $\C$ admits an on-site symmetry realization on $\mathfrak A_q$ for some positive integer $q$ if and only if it is tensor-generated by a connected Frobenius algebra of depth at most two. The algebra may be chosen to be a Q-system. If $\C\not\simeq\Vect_{\mathbb C}$, its depth is exactly two.
\end{corollary}
\end{samepage}
\begin{proof}
By \cite[Corollary~4.5]{EJ}, on-site realizability is equivalent to the existence of a unitary fiber functor. An ordinary fiber functor makes $\Vect_{\mathbb C}$ a semisimple module category with one simple object. This module category is unitarizable by \cite[Corollary~3.8]{CGGH}; its action then gives a unitary fiber functor to finite-dimensional Hilbert spaces. Conversely, forgetting inner products gives an ordinary fiber functor. Apply Corollary~\ref{cor:frobenius}. A connected Frobenius algebra is algebra-isomorphic to a Q-system by \cite[Theorem~3.2 and Remark~3.3]{CGGH}; algebra isomorphisms preserve generation and depth. The final assertion is Corollary~\ref{cor:exact-depth}.
\end{proof}

The criterion expresses on-site realizability through one algebra object and the module identity $A\otimes A\cong A^{\oplus n}$. It is stronger than integrality, which characterizes the more general tensor-product realizations in \cite[Proposition~4.1 and Theorem~4.2]{EJ}. It also gives an obstruction under quantum cellular automata, that is, $*$-automorphisms for which both the map and its inverse enlarge supports by a uniformly bounded amount. Such a change of variables transports the symmetry category by a tensor equivalence \cite[Remark~3.19]{EJ}. If the category has no connected tensor-generating algebra of depth at most two, no transported realization can be on-site. The criterion concerns existence; it does not assert that every realization of a category admitting a fiber functor is related to an on-site one by a quantum cellular automaton.

\subsection{A nonsemisimple example without normalization}
In the following nonsemisimple example of Corollary~\ref{cor:frobenius}, the regular comodule algebra has depth two even though its Frobenius composites $\lambda\circ e$ and $m\circ d$ both vanish and hence cannot be normalized to the identity.

\begin{example}\label{ex:charp}
Let $\operatorname{char}(k)=p>0$ and let
\[
H=k[t]/(t^p),\qquad
\Delta_H(t)=t\otimes1+1\otimes t,\quad
\eps_H(t)=0,\quad S(t)=-t.
\]
The identity $(t\otimes1+1\otimes t)^p=0$ makes this a Hopf algebra. Define $\lambda(t^r)=\delta_{r,p-1}$ for $0\le r<p$. Expanding $\Delta_H(t^r)$ shows that
$((\id_H\otimes\lambda)\circ\Delta_H)(h)=\lambda(h)1$. The pairing $\lambda(ab)$ has an invertible antidiagonal matrix in the basis $1,t,\ldots,t^{p-1}$. Its inverse copairing and Frobenius coproduct are
\[
c=\sum_{r=0}^{p-1}t^r\otimes t^{p-1-r},
\qquad
d(t^m)=\sum_{r=0}^{p-1-m}t^{m+r}\otimes t^{p-1-r}.
\]
Nevertheless $\lambda(1)=0$ and $m_H\circ d=0$. For the latter assertion, when $m>0$ every product in the displayed sum is zero, and when $m=0$ the result is $p\,t^{p-1}=0$. Thus the Frobenius structure cannot be normalized by either of these scalars. The regular comodule still satisfies
$H\otimes H\cong H^{\oplus p}$ as a right $H$-module in $\Corep(H)$, by Proposition~\ref{prop:hopf-depth}.

This category is not semisimple and has only one simple object. Indeed, if $f_r\in H^*$ extracts the coefficient of $t^r$, convolution gives $f_1^r=r!f_r$ for $0\le r<p$ and $f_1^p=0$. Hence $H^*\cong k[u]/(u^p)$ as an algebra. In particular the category is not equivalent to $\Vect_k$. Propositions~\ref{prop:hopf-depth} and \ref{prop:depth-two-identity}, together with Corollary~\ref{cor:exact-depth}, give $\depth(H)=2$. The example shows that neither $\lambda\circ e\ne0$ nor $m\circ d\ne0$ is needed in the depth-two characterization.
\end{example}

\section*{Acknowledgments}
Mainak Ghosh is supported by a Start-up grant from Hetao Institute of
Mathematics and Interdisciplinary Sciences.
S\'ebastien Palcoux is supported by the National Natural Science
Foundation of China (NSFC), General Program (Grant no.~12471031).
GPT-6 Astra served as a research assistant; the authors take full
responsibility for the final manuscript.

\end{document}